\documentclass{conm-p-l}

\usepackage{eucal}%

\usepackage{hyperref}%

\usepackage{tikz-cd}%

\usepackage{thmtools}

\declaretheorem[style=theorem,qed=\qedsymbol,numberwithin=section]{theorem}
\declaretheorem[style=theorem,qed=\qedsymbol,sibling=theorem]{proposition}
\declaretheorem[style=definition,qed=\qedsymbol,sibling=theorem]{definition}
\declaretheorem[style=remark,qed=\qedsymbol,sibling=theorem]{remark}
\declaretheorem[style=remark,qed=\qedsymbol,sibling=theorem]{example}

\usepackage{xparse}

\NewDocumentCommand{\morphism}{o m m o O{normal}}{%
  \begin{tikzcd}[ampersand replacement=\&, column sep=#5]
    \IfValueT{#1}{#1\colon}#2\IfValueTF{#4}{\rar{#4}}{\rar}\&#3
  \end{tikzcd}
}

\newcommand{\MyCategory}[1]{\mathbf{#1}}

\NewDocumentCommand{\SimplexCat}{o o}{\MyCategory{\Delta}\IfValueT{#1}{(#1,#2)}}

\NewDocumentCommand{\ord}{s m}{[#2]\IfBooleanT{#1}{_+}}

\newcommand{\coloneqq}{:=}

\NewDocumentCommand{\set}{m o}{\left\{\,#1%
    \IfNoValueF{#2}{\,\middle|\,#2}%
    \,\right\}%
}

\renewcommand{\AA}{\mathbb{A}}

\newcommand{\zahlen}[1]{\mathbb{#1}}
\newcommand{\ZZ}{\zahlen{Z}}

\newcommand{\category}[1]{\mathcal{#1}}
\newcommand{\A}{\category{A}}
\newcommand{\C}{\category{C}}

\newcommand{\op}{{\mathrm{op}}}

\NewDocumentCommand{\adjunction}{o m m o O{} O{normal}}{
  \begin{tikzcd}[ampersand replacement=\&, column sep=#6]
     \IfValueT{#1}{#1\colon} #2\rar[shift left]\& #3\IfValueT{#4}{\noloc #4}\lar[shift left]
  \end{tikzcd}
}

\newcommand{\AbelianGroups}{\MyCategory{Ab}}
\newcommand{\sAb}{\AbelianGroups_{\SimplexCat}}

\newcommand{\noloc}{\,\colon\!\!}

\NewDocumentCommand{\GeoRealisation}{O{-} o}{%
  |#1|\IfValueT{#2}{_{#2}}
}

\newcommand{\Sets}{\MyCategory{Set}}

\newcommand{\sSets}{\Sets_{\SimplexCat}}

\NewDocumentCommand{\Sing}{s o}{%
  \operatorname{Sing}\IfValueT{#2}{\IfBooleanTF{#1}{#2}{\left(#2\right)}}
}

\NewDocumentCommand{\Hom}{O{-} O{-} o}{%
  \operatorname{Hom}\IfValueT{#3}{_{#3}}\left({#1},{#2}\right)%
}

\newcommand{\Stovicek}{{{\v{S}}\v{t}ov\'{i}\v{c}ek}}

\NewDocumentCommand{\Map}{O{-} O{-} o}{%
  \operatorname{Map}\IfValueT{#3}{_{#3}}({#1},{#2})%
}

\NewDocumentCommand{\Ho}{s o}{%
  \operatorname{Ho}\IfValueT{#2}{\IfBooleanTF{#1}{#2}{\left(#2\right)}}
}

\NewDocumentCommand{\Fun}{o o o o}{%
  \operatorname{Fun}\IfValueT{#3}{^{#3}}\IfValueT{#4}{_{#4}}\IfValueT{#1}{(#1,#2)}%
}

\NewDocumentCommand{\functor}{s o m m o O{normal}}{%
  \begin{tikzcd}[ampersand replacement=\&, column sep=#6]
    \IfValueT{#2}{#2\colon} #3\IfBooleanTF{#1}{\IfValueTF{#5}{\rar[mapsto]{#5}}{\rar[mapsto]}}{\IfValueTF{#5}{\rar{#5}}{\rar}}\& #4
  \end{tikzcd}
}

\newcommand{\sS}{\mathsf{S}}%
\newcommand{\Sm}{\sS^{\langle m\rangle}}%

\newcommand{\Ab}{\AbelianGroups}%

\newcommand{\setP}[2]{\set{#1}[#2]}

\newcommand{\recollementnoph}[6]{%
  \begin{tikzcd}[ampersand replacement=\&, column sep=large]%
    #1 \ar{r}[description]{#5} \&%
    #2 \ar{r}\ar[shift right=4]{l}[description]{#4} \ar[shift left=4]{l}[description]{#6}\&%
    #3 \ar[shift right=4]{l} \ar[shift left=4]{l}%
  \end{tikzcd}%
}

\title[Higher Auslander algebras of type $\AA$ and higher
$\sS$-cons\-truct\-ions]{Higher Auslander algebras of type $\AA$\\ and the higher
  Waldhausen $\sS$-cons\-truct\-ions}

\author[G. Jasso]{Gustavo Jasso}

\address{Mathematisches Institut\\
Rheinische Friedrich-Wilhelms-Universität Bonn\\
Endenicher Allee 60\\
D-53115 Bonn, Germany}

\email{gjasso@math.uni-bonn.de}
\urladdr{http://gustavo.jasso.info}

\keywords{Auslander--Reiten theory; Eilenberg--Mac Lane spaces; Waldhausen
  $\sS$-construction}

\subjclass[2010]{Primary: 18G30. Secondary 16G70, 18E30}

\begin{document}

\maketitle

\begin{abstract}
  These notes are an expanded version of my talk at the ICRA 2018 in Prague,
  Czech Republic; they are based on joint work with Tobias Dyckerhoff and Tashi
  Walde \cite{DJW19b}. In them we relate Iyama's higher Auslander algebras of
  type $\AA$ to Eilenberg--Mac Lane spaces in algebraic topology and to
  higher-di\-men\-sio\-nal versions of the Waldhausen $\sS$-cons\-truct\-ion
  from algebraic $K$-theory.
\end{abstract}

\section{Motivation: Eilenberg--Mac Lane spaces}

Let $A$ be an abelian group. For each positive integer $m$ there exists a
topological space $K(A,m)$ characterised, up to weak homotopy
equivalence\footnote{A \emph{weak homotopy equivalence} is a continuous map
  $X\to Y$ between topological spaces which induces a bijection
  $$\pi_0(X)\xrightarrow{\cong}\pi_0(Y)$$ on connected components as well as group isomorphisms
  $$\pi_n(X,x)\xrightarrow{\cong}\pi_n(Y,f(x))$$ for all points $x\in X$ and all $n\geq 1$.
  Two topological spaces are \emph{weakly homotopy equivalent} if they are
  connected by a zig-zag of weak homotopy equivalences.}, by the existence of
natural bijections\footnote{Here $[X,K(A,m)]$ denotes the set of homotopy
  classes of continuous maps $X\to K(A,m)$, while $H^m(X;A)$ denotes the $m$-th
  singular cohomology group of $X$ with coefficients in $A$.}
\begin{equation}
  \label{eq:Kam}
  \morphism{\left[X,K(A,m)\right]}{H^m(X;A)}[\cong]
\end{equation}
where $X$ is a CW complex, see Corollary III.2.17 in \cite{GJ99}. Equivalently,
the topological space $K(A,m)$ is characterised by the following two properties:
\begin{itemize}
\item The topological space $K(A,m)$ is path connected;
\item there is an isomorphism
  \begin{equation*}
    \pi_n(K(A,m))\cong\begin{cases}
      A&\text{if }n=m,\\
      0&\text{otherwise}.
    \end{cases}
  \end{equation*}
\end{itemize}
A topological space satisfying the above properties is called an
\emph{Ei\-len\-berg--Mac Lane space}\footnote{The symbol $K(A,m)$ is more
  commonly used to denote the weak homotopy equivalence class (also known as the
  homotopy type) of an Eilenberg--Mac Lane space. For our purposes it will be
  convenient to abuse the notation and use this symbol to denote a preferred
  representative of this equivalence class.}. The bijection \eqref{eq:Kam}
exhibits Eilenberg--Mac Lane spaces as fundamental objects in algebraic topology
while their characterisation in terms of higher homotopy groups makes their
relative simplicity manifest.

Rather surprisingly, the standard construction of the Eilenberg--Mac Lane
spac\-es can be expressed in terms of the \emph{higher Auslander--Reiten theory}
of Iyama's higher Auslander algebras of type $\AA$. This observation, which we
make precise below, is one of the starting points of our investigations.

\subsection{Interlude: the simplex category}

The modern approach to homotopy theory makes extensive use of simplicial methods
and the standard construction of Ei\-len\-berg--Mac Lane spaces is no exception.
We recall a minimal amount of terminology from simplicial homotopy
theory\footnote{The reader is referred to \cite{GJ99} for an introduction to
  simplicial homotopy theory.} needed for describing this construction.

We remind the reader of the definition of the \emph{simplex category}, which is
commonly denoted by $\SimplexCat$. The objects of the simplex category are the
linear posets
\begin{equation*}
  \ord n\coloneqq\set{0\to 1\to\cdots\to n},
\end{equation*}
where $n$ ranges over all non-negative integers, with morphisms the monotone
functions $\ord m\to\ord n$. We denote the set of morphisms $\ord m\to\ord n$ in
$\SimplexCat$ by $\SimplexCat(m,n)$. Note that the set $\SimplexCat(m,n)$ has a
natural poset structure\footnote{The poset structure on the morphism sets in
  $\SimplexCat$ allows us to promote the simplex category to a $2$-category;
  this perspective, implicit in parts of these notes, is essential in our
  work \cite{DJW19b}.}: Namely, $\sigma\leq \tau$ if and only if $\sigma_i\leq\tau_i$ for all
$i\in\ord m$. The elements of $\SimplexCat(m,n)$ are called \emph{$m$-simplices
  in $\SimplexCat^n$}; an $m$-simplex is \emph{non-degenerate} if the underlying
monotone function is injective and is \emph{degenerate} otherwise. We often
identify $\SimplexCat(m,n)$ with the set of ordered tuples $\sigma\in\ZZ^{m+1}$
which satisfy the inequalities $0\leq
\sigma_0\leq\sigma_1\leq\cdots\leq\sigma_m\leq n$.

A \emph{simplicial object}\footnote{Depending on the target category one speaks
  of simplicial sets, simplicial topological spaces, simplicial abelian groups,
  etc.} in a category $\C$ is a functor $X\colon\SimplexCat^\op\to\C$. It is
elementary to verify, using a specific presentation of $\SimplexCat$ in terms of
generators and relations, that the data of a simplicial object in $\C$ is
equivalent to that of a diagram in $\C$ of the form
\begin{equation*}
  \begin{tikzcd}
    \cdots\rar[shift right=4.5]\rar[shift left=4.5]\rar[shift
    right=1.5]\rar[shift left=1.5]&X_2\rar\rar[shift right=3]\rar[shift
    left=3]\lar\lar[shift right=3]\lar[shift left=3]&X_1\rar[shift
    left=1.5]\rar[shift right=1.5]\lar[shift left=1.5]\lar[shift
    right=1.5]&X_0\lar
  \end{tikzcd}
\end{equation*}
where $X_n\coloneqq X(\ord n)$ is the object of \emph{$n$-simplices of $X$}. The
morphisms
\begin{equation*}
  \setP{d_i\colon X_n \to X_{n-1}}{i\in\ord n},
\end{equation*}
called \emph{face maps}, and the morphisms
\begin{equation*}
  \setP{s_i\colon X_n\to X_{n+1}}{i\in\ord n},
\end{equation*}
called \emph{degeneracy maps}, are subject to a number of relations commonly
referred to as the \emph{simplicial identities}\footnote{Let $X$ be a simplicial
  set. The simplicial identities allow us to think of an element of $X_n$ as an
  abstract (possibly degenerate) $n$-simplex whose boundary consists of the
  $(n-1)$-simplices $\setP{d_i(x)\in X_{n-1}}{i\in[n]}$.}:
\begin{equation*}
  \begin{aligned}[c]
    d_i\circ d_j&=d_{j-1}\circ d_i\quad i<j,\\
    s_i\circ s_j&=s_j\circ s_{i-1}\quad i>j,\\
  \end{aligned}\qquad\quad
  \begin{aligned}[c]
    d_i\circ s_j&=\begin{cases}
      s_{j-1}\circ d_i&i<j,\\
      s_j\circ d_{i-1}&i>j+1,\\
      1&\text{otherwise}.
    \end{cases}
  \end{aligned}
\end{equation*}
Thus, a simplicial object in $\C$ is a \emph{contravariant representation} in
the category $\C$ of the infinite quiver with (inadmissible) relations
corresponding to the standard presentation of $\SimplexCat$.

\subsection{Construction via the Dold--Kan correspondence}

The \emph{Dold--Kan correspondence} \cite{Dol58,Kan58} is an explicit adjoint
equivalence of categories\footnote{A detailed discussion of this equivalence can
  be found in Section III.2 in \cite{GJ99}.}
\begin{equation*}
  \adjunction[C]{\sAb}{\operatorname{Ch}_{\geq0}(\Ab)}[N]
\end{equation*}
between the (abelian) category of simplicial abelian groups and the (abelian)
category of chain complexes of abelian groups which are concentrated in
non-negative homological degrees\footnote{From a representation-theoretic
  perspective, the Dold--Kan correspondence should be viewed as a Morita
  equivalence between $\ZZ\SimplexCat$, the $\ZZ$-linear envelope of
  $\SimplexCat$, and the $\ZZ$-linear path category of the quiver $0\to 1\to
  2\to \cdots$ modulo the square of the Jacobson radical (note that we work with
  \emph{contravariant} modules).}.

Let $A$ be an abelian group and $m$ a positive integer. As usual, let us denote
by $A[m]$ the chain complex of abelian groups whose only non-zero component
consists of $A$ placed in homological degree $m$. The Eilenberg--Mac Lane space
$K(A,m)$ is defined as
\begin{equation*}
  K(A,m)\coloneqq|N(A[m])_\bullet|,
\end{equation*}
that is as the geometric realisation\footnote{The geometric realisation functor
  is part of an adjunction \cite[p.~9]{GJ99}
  $$\adjunction[\GeoRealisation]{\sSets}{\mathbf{CGHaus}}[\Sing]$$
  between the category of simplicial sets and the category of compactly
  generated Hausdorff spaces. In fact, this adjunction is a Quillen equivalence
  with respect to appropriate model category structures \cite{Qui67,Hov99}. This
  Quillen equivalence justifies the simplicial approach to homotopy theory.} of
the underlying simplicial set of $N(A[m])_\bullet$.

\subsection{The Eilenberg--Mac Lane space $K(A,1)$}

Let $A$ be an abelian group. For our purposes it is instructive to give an
explicit description of the simplicial abelian group $N(A[m])_\bullet$ in the
simplest case $m=1$.

\begin{definition}
  The simplicial abelian group $N(A[1])_\bullet$ is defined as follows. For
  $n\geq0$, let $N(A[1])_n$ be the abelian group of upper-triangular arrays
  \begin{equation*}
    \label{eq:element_NA1n}
    \begin{pmatrix}
      a_{00}&a_{01}&a_{02}&\cdots&a_{0,n-1}&a_{0n}\\
      &a_{11}&a_{12}&\cdots&a_{1,n-1}&a_{1n}\\
      &&\ddots&\ddots&\vdots&\vdots\\
      &&&&a_{n-1,n-1}&a_{n-1,n}\\
      &&&&& a_{nn}
    \end{pmatrix}
  \end{equation*}
  with entries in $A$ such that for each $0\leq i\leq n$ there is an equality
  $a_{ii}=0$ and for all $0\leq i<j<k\leq n$ the \emph{Euler relation}
  \begin{equation}
    \label{eq:Euler_relations_NA1}
    a_{ij} - a_{ik} + a_{jk}=0
  \end{equation}
  is satisfied. The $i$-th face map
  \[
    d_i\colon N(A[1])_n\to N(A[1])_{n-1}
  \]
  is given by deleting the $i$-th row and the $i$-th column while the $i$-th
  degeneracy map
  \[
    s_i\colon N(A[1])_n\to N(A[1])_{n+1}
  \]
  is given by repeating the $i$-th row and the $i$-th column.
\end{definition}
  
\begin{remark}
  An $n$-simplex of $N(A[1])_\bullet$ is completely determined by the $n$-tuple
  $(a_{01},a_{02},\dots,a_{0n})$. Indeed, the Euler relations
  \eqref{eq:Euler_relations_NA1} imply that for each $1\leq j<k\leq n$ the
  equality
  \begin{equation*}
    a_{jk}=a_{0k}-a_{0j}
  \end{equation*}
  is satisfied. We conclude that there is an isomorphism
  \begin{equation*}
    N(A[1])_n\cong A^n.
  \end{equation*}
  Note, however, that the above identification makes the simplicial structure of
  $N(A[1])_\bullet$ less apparent.
\end{remark}

Our next task is to relate the previous construction of the Ei\-len\-berg--Mac
Lane space $K(A,1)$ to the Auslander--Reiten theory of the family of quivers
\begin{align*}
  \AA_n&\coloneqq1\to2\to\cdots\to n;&&n\geq0.
\end{align*}
Consider the Auslander--Reiten quiver of the category of
finite-di\-men\-sio\-nal representations\footnote{We work over some fixed but
  unspecified field.} of the quiver $\AA_n$, which we depict as follows:
\begin{equation*}
  \begin{tikzcd}[column sep=tiny,row sep=tiny]
    M_{00}\rar&M_{01}\rar\dar&M_{02}\rar\dar&\cdots\rar&M_{0,n-1}\rar\dar&M_{0n}\dar\\
    &M_{11}\rar&M_{12}\rar\dar&\cdots\rar&M_{1,n-1}\rar\dar&M_{1n}\dar\\
    &&\ddots&\ddots&\vdots\dar&\vdots\dar\\
    &&&&M_{n-1,n-1}\rar&M_{n-1,n}\dar\\
    &&&&& M_{nn}
  \end{tikzcd}
\end{equation*}
For reasons which shall become clear shortly, we have included additional
vertices $M_{ii}=0$ in the above Auslander--Reiten quiver. The Grothendieck
group $K_0(\AA_n)$ admits a presentation by generators and relations as the
quotient of the free abelian group with basis
\begin{equation*}
  \setP{[M_{ij}]}{0\leq i\leq j\leq n}
\end{equation*}
modulo the subgroup generated by the relations $[M_{ii}]=0$, $i\in\ord n$,
together with the Euler relations
\begin{equation}
  \label{eq:Euler_relations_K0An}
  [M_{ij}]-[M_{ik}]+[M_{jk}]=0
\end{equation}
corresponding to short exact sequences\footnote{Equivalently, the Euler
  relations are induced by biCartesian squares
  $$
  \begin{tikzcd}[ampersand replacement=\&]
    M_{ij}\rar\dar\&M_{ik}\dar\\
    M_{ii}\rar\&M_{jk}
  \end{tikzcd}
  $$
  where $0\leq i<j<k\leq n$.}
\begin{equation*}
  0\to M_{ij}\to M_{ik}\to M_{jk}\to0
\end{equation*}
in the abelian category of finite-dimensional representations of the quiver
$\AA_n$, where $0\leq i<j<k\leq n$.

It is immediate from the above presentation of $K_0(\AA_n)$ that there are
isomorphisms
\begin{equation}
  \label{eq:Na1n}
  N(A[1])_n\cong\Hom[K_0(\AA_n)][A][\ZZ]\cong\Hom[\ZZ^n][A][\ZZ]\cong A^n,
\end{equation}
where $n\geq 1$. We leave it to the reader to verify that the Grothendieck
groups $K_0(\AA_n)$ assemble into a co-simplicial abelian group, that is into a
functor
\begin{equation*}
  K_0(\AA_\bullet)\colon\SimplexCat\to\Ab,
\end{equation*}
where $K_0(\AA_0)\coloneqq0$. Moreover, the isomorphisms \eqref{eq:Na1n}
assemble into an isomorphism of simplicial abelian groups
\begin{equation}
  \label{eq:NA1_corepresentable}
  N(A[1])_\bullet\cong\Hom[K_0(\AA_\bullet)][A][\ZZ]
\end{equation}
where $\Hom[K_0(\AA_\bullet)][A][\ZZ]$ denotes the composite
\begin{equation*}
  \begin{tikzcd}[ampersand replacement=\&, column sep=5em]   
    \SimplexCat^\op\rar{K_0(\AA_\bullet)^\op}\&\Ab^\op\rar{\Hom[-][A][\ZZ]}\&\Ab,
  \end{tikzcd}
\end{equation*}
and $K_0(\AA_\bullet)^\op$ indicates the passage to the opposite functor.

\begin{remark}
  It is elementary to verify that the Euler relations
  \eqref{eq:Euler_relations_K0An} are generated by the \emph{Auslander--Reiten
    relations}
  \begin{equation*}
    [M_{ij}]-([M_{i+1,j}]+[M_{i,j+1}])+[M_{i+1,j+1}]=0
  \end{equation*}
  corresponding to the almost-split sequences\footnote{Equivalently, the
    Auslander--Reiten relations are induced by meshes
    $$
    \begin{tikzcd}[ampersand replacement=\&]
      M_{ij}\rar\dar\&M_{i,j+1}\dar\\
      M_{i+1,j}\rar\&M_{i+1,j+1}
    \end{tikzcd}
    $$
    in the Auslander--Reiten quiver, where $0\leq i<j<n$.}
  \begin{equation*}
    0\to M_{ij}\to M_{i+1,j}\oplus M_{i,j+1}\to M_{i+1,j+1}\to0,
  \end{equation*}
  where $0\leq i<j<n$, see also \cite{Aus84}. In particular, the isomorphism
  \eqref{eq:NA1_corepresentable} gives a precise relationship between the
  simplicial abelian group $N(A[1])_\bullet$ and the Auslander--Reiten theory of
  the linearly oriented quivers of Dynkin type $\AA$.
\end{remark}

\subsection{The Eilenberg--Mac Lane space $K(A,m)$}
  
Let $A$ be an abelian group and $m$ a positive integer. The previous discussion
extends to the Eilenberg--Mac Lane space $K(A,m)$ provided that we replace the
linearly oriented quivers of type $\AA$ by Iyama's $m$-dimensional Auslander
algebras $\AA_\ell^{(m)}$ of type $\AA$, see \cite{Iya11} for the definition.
Indeed, the Grothendieck groups of these algebras assemble into a co-simplicial
abelian group
\begin{equation*}
  K_0(\AA_{\bullet-m+1}^{(m)})\colon\SimplexCat\to\Ab,
\end{equation*}
where $K_0(\AA_{n-m+1}^{(m)})\coloneqq0$ for $n<m$. In analogy with the case
$m=1$, there is a canonical isomorphism of simplicial abelian groups
\begin{equation*}
  N(A[m])_\bullet\cong\Hom[K_0(\AA_{\bullet-m+1}^{(m)})][A][\ZZ],
\end{equation*}
see Theorem 1.20 in \cite{DJW19b} for a detailed description of the above
isomorphism.

The above description of $N(A[m])_\bullet$ utilises the presentation of the
Grothendieck group $K_0(\AA_{n-m+1}^{(m)})$ in terms of the basis
\begin{equation*}
  \setP{[M_\sigma]}{\sigma\in\SimplexCat(m,n)}
\end{equation*}
labelled by the indecomposable summands of the unique basic $m$-cluster-tilting
$\AA_n^{(m)}$-module\footnote{Implicit in this discussion is the Grothendieck
  group of an $m$-cluster-tilting subcategory, which can be defined using the
  ideas in \cite{BT14}.}. In analogy with the case $m=1$, one imposes relations
$[M_\sigma]=0$ for all degenerate $m$-simplices in $\SimplexCat^n$ as well as
\emph{higher Euler relations}
\begin{equation}
  \label{eq:Euler_relations_Am_1}
  \sum_{i=0}^m(-1)^i[M_{d_i(\sigma)}]=0,
\end{equation}
where $\sigma$ is a non-degenerate $(m+1)$-simplex in $\SimplexCat^n$ and
\begin{equation*}
  d_i(\sigma)=(\sigma_0,\sigma_1,\cdots,\sigma_{i-1},\widehat{\sigma}_i,\sigma_{i+1},\cdots,\sigma_m).
\end{equation*}
The higher Euler relations are induced by exact sequences
\[
  0\to M_{d_m(\sigma)}\to M_{d_{m-1}(\sigma)}\to\cdots\to M_{d_1(\sigma)}\to
  M_{d_0(\sigma)}\to0
\]
in the abelian category of finite-dimensional representations of the higher
Auslander algebra $\AA_{n-m+1}^{(m)}$. Finally, we note that the higher Euler
relations are generated by the \emph{higher Auslander--Reiten relations}
\begin{equation}
  \label{eq:Euler_relations_Am}
  \sum_{v}(-1)^{|v|}[M_{\sigma+v}]=0,
\end{equation}
where $\sigma$ is a non-degenerate $m$-simplex in $\SimplexCat^n$ with
$\sigma_m<n$, the tuple $v$ ranges over all vertices of the $(m+1)$-cube
\begin{equation*}
  \underbrace{\set{0,1}\times\cdots\times\set{0,1}}_{m+1\text{ times}},
\end{equation*}
and $|v|\coloneqq v_0+v_1+\cdots+v_m$. Again in analogy with the case $m=1$, the
relations \eqref{eq:Euler_relations_Am} are induced by the $m$-almost-split
sequences
\[
  0\to
  M_\sigma\to\bigoplus_{|v|=1}M_{\sigma+v}\to\cdots\to\bigoplus_{|v|=m}M_{\sigma+v}\to
  M_{\sigma+(1,\dots,1)}\to0
\]
in the unique $m$-cluster-tilting subcategory associated with
$\AA_{n-m+1}^{(m)}$, see \cite{Iya07,Iya11} and \cite{OT12} for definitions and
further details.

\subsection{On the purpose of these notes}

In the remainder of these notes is to describe a \emph{categorified} version of
the above discussion where the abelian group $A$ is replaced by a refined
version of a triangulated category\footnote{Triangulated categories are
  inadequate for this purpose as the cone of a morphisms lacks the necessary
  functoriality.}, namely a \emph{stable $\infty$-category} in the sense of
Lurie \cite{Lur17}. As such, the simplicial objects we shall consider can be
thought of as `categorified Eilenberg--Mac Lane spaces' as originally advocated
by Dyckerhoff in \cite{Dyc17}. The results presented in these notes can be also
seen as contributions to the abstract representation theory program of Groth and
\Stovicek~\cite{GS18b}, see also Ladkani's Ph.D.~Thesis \cite{Lad08}. Indeed,
the results presented here are for the most part higher-dimensional versions of
results in \cite{GS16}.

\section{Rudiments of the theory of $\infty$-categories}

The language of $\infty$-categories, developed by Joyal in \cite{Joy02} as well
as in unpublished work and by Lurie in \cite{Lur09,Lur17}, affords an adequate
framework for our results. It is however impractical to include here a formal
introduction to the theory of $\infty$-categories\footnote{The standard
  reference for the theory of $\infty$-categories is Lurie's book \cite{Lur09}.
  Alternatives include Groth's lecture notes \cite{Gro10} and Cisinski's book
  \cite{Cis19}.}. Instead we aim to provide the reader with a minimal amount of
intuition regarding the theory which we hope will be sufficient to follow our
exposition in the subsequent section.

\subsection{What are $\infty$-categories?}

An $\infty$-category is a mathematical structure\footnote{More precisely, an
  $\infty$-category is a particular kind of simplicial set called a `weak Kan
  complex'.} which implements the idea of a
\begin{quote}
  `higher-dimensional category with morphisms of every positive degree endowed
  with a coherently-associative composition law and such that the morphisms of
  degree greater than 1 are invertible'.
\end{quote}
In other words the theory of $\infty$-categories is a model for the theory of
$(\infty,1)$-categories. In particular, for every pair of objects $x$ and $y$ in
an $\infty$-category $\C$ there is a `space' $\Map[x][y][\C]$ of maps $x\to y$.
Within this paradigm ordinary categories are those for which the non-empty
spaces of maps are homotopy equivalent to discrete spaces\footnote{See
  Propositions 2.3.4.5 and 2.3.4.18 in \cite{Lur09}.}. After suitable
translation, further examples of $\infty$-categories are provided by
differential graded categories and, more generally,
$A_\infty$-categories\footnote{See Section 1.3 in \cite{Lur17} for the case of
  differential graded categories and \cite{Fao17a,Fao17b} for the case of
  $A_\infty$-categories.}.

An important feature of the theory of $\infty$-categories is that it allows us
to formalise the notion of a universal property which only holds `up to coherent
homotopy'. For example, an object $x$ of an $\infty$-category $\C$ is
\emph{initial} if for each object $y$ of $\C$ the space $\Map[x][y][\C]$ is
contractible (whence in particular non-empty).

Just as in ordinary category theory, universal properties are also captured by
appropriate notions of (homotopy) limit and (homotopy) colimit\footnote{In
  $\infty$-category theory it is customary to drop the qualifier `homotopy' from
  the terminology as only homotopy-invariant notions make sense in this
  context.}. For example, given two morphisms $x\to y'$ and $x\to y''$ in an
$\infty$-category $\C$, we may define the (homotopy) pushout of the
diagram\footnote{Strictly speaking, specifying a coherently commutative diagram
  in an $\infty$-category involves an infinite amount of data. For expository
  purposes we nonetheless display such diagrams as we would display commutative
  diagrams in ordinary categories.}
\begin{equation*}
  \begin{tikzcd}
    x\rar\dar&y'\\
    y''
  \end{tikzcd}
\end{equation*}
by means of a suitable universal property, that is as a (homotopy) colimit
diagram of the form\footnote{We decorate the square to indicate that it is a
  (homotopy) pushout square.}
\begin{equation*}
  \begin{tikzcd}
    x\rar\dar\ar[phantom]{dr}[description,near end]{\ulcorner}&y'\dar\\
    y''\rar&z
  \end{tikzcd}
\end{equation*}
The space of all possible (homotopy) (co)limit diagrams of a fixed shape is
either empty (if no (co)limit of the diagram exists in $\C$) or
contractible\footnote{This means that, for the purposes of $\infty$-category
  theory, (co)limit cones of a fixed diagram in an $\infty$-category are unique
  in the appropriate sense. See Proposition 1.2.12.9 and Definition 1.2.13.4 in
  \cite{Lur09} for details.}.

\subsection{Stable $\infty$-categories}

In the sequel we are mostly concerned with the following class of
$\infty$-categories, which the reader might want to think of as refined versions
of triangulated categories.

\begin{definition}
  An $\infty$-category $\A$ is \emph{stable}\footnote{Although we have not
    provided formal definitions of any the notions involved, the property of
    being \emph{stable} is sufficiently intuitive for it to be worth to be
    included in these notes. A detailed treatment of the theory of stable
    $\infty$-categories can be found in Chapter 1 in \cite{Lur17}.} if it has
  the following properties:
  \begin{enumerate}
  \item The $\infty$-category $\A$ is \emph{pointed}, that is $\A$ has a zero
    object.
  \item For every morphism $f\colon x\to y$ in $\A$ there exist (coherently
    commutative) squares in $\A$ of the form
    \begin{equation*}
      \begin{tikzcd}
        x\rar{f}\dar\ar[phantom]{dr}[description,near end]{\ulcorner}&y\dar\\
        0\rar&z
      \end{tikzcd}\qquad{and}\qquad
      \begin{tikzcd}
        w\rar\dar\ar[phantom]{dr}[description,near start]{\lrcorner}&x\dar{f}\\
        0\rar&y
      \end{tikzcd}
    \end{equation*}
    which are a (homotopy) pushout square and a (homotopy) pullback square,
    respectively. The objects $z$ and $w$ are called the \emph{cofibre of $f$} and
    the \emph{fibre of $f$}, respectively\footnote{The cofibre and the fibre of
      a morphism in a stable $\infty$-category are the $\infty$-categorical
      versions of the cone and the co-cone of a morphism in a triangulated
      category. Following the established convention, we adopt the `topological'
      terminology which further reminds us of the homotopical nature of these
      concepts.}.
  \item A diagram in $\A$ of the form
    \begin{equation*}
      \begin{tikzcd}
        x\rar\dar&y\dar\\
        0\rar&z
      \end{tikzcd}
    \end{equation*}
    is a (homotopy) pushout square if and only if it is a (homotopy) pullback
    square. A diagram as above which satisfies these additional properties is
    called a \emph{fibre-cofibre sequence}\footnote{Fibre-cofibre sequences in stable
      $\infty$-categories are the $\infty$-categorical versions of exact
      triangles in triangulated categories.}.\qedhere
  \end{enumerate}
\end{definition}

\begin{remark}
  Stable $\infty$-categories are defined in terms of \emph{properties}
  and not in terms of additional \emph{structure}, the latter being the case for
  triangulated categories. Roughly speaking, this is the reason why stable
  $\infty$-categories enjoy better \emph{formal properties} than triangulated
  categories.
\end{remark}

Every $\infty$-category $\C$ has an associated \emph{homotopy category}
$\Ho[\C]$ which is in fact an ordinary category. The following result\footnote{\autoref{prop:htpy_cat_triangulated} should be compared with Happel's theorem
  which states that the stable category of a Frobenius exact category is a
  triangulated category \cite{Hap88}.}
relates stable $\infty$-categories to triangulated categories, see Theorem
1.1.2.14 in \cite{Lur17}.
\begin{proposition}
  \label{prop:htpy_cat_triangulated} Let $\A$ be a stable $\infty$-category.
  Then, the homotopy category $\Ho[\A]$ is (canonically) a triangulated
  category.
\end{proposition}

\begin{remark}
  The suspension $\Sigma(x)$ of an object $x$ of a stable $\infty$-category $\A$
  is characterised by the existence of a fibre-cofibre sequence\footnote{The
    decoration indicates that the square is both a (homotopy) pushout and a
    (homotopy) pullback. In fact, in a stable $\infty$-category these two
    notions coincide, see Proposition 1.1.3.4 in \cite{Lur17}.}
    $$
    \begin{tikzcd}[ampersand replacement=\&]
      x\rar\dar\ar[phantom]{dr}[description]{\square}\&0\dar\\
      0\rar\&\Sigma(x)
    \end{tikzcd}
    $$
    Similarly, exact triangles in the homotopy category $\Ho[\A]$ are induced
    by diagrams in $\A$ of the form
    $$
    \begin{tikzcd}[ampersand replacement=\&]
      x\rar\dar\ar[phantom]{dr}[description]{\square}\&y\dar\rar\ar[phantom]{dr}[description]{\square}\&0\dar\\
      0\rar\&z\rar\&\Sigma(x)
    \end{tikzcd}
    $$
    in which each square is a fibre-cofibre sequence.
  \end{remark}

  The good formal behaviour of stable $\infty$-categories is illustrated
  by the following statement\footnote{ \autoref{prop:pointwise} should be
    compared with the elementary fact that the category of functors from a small
    category into an abelian category is again abelian. Indeed, (homotopy)
    limits and (homotopy) colimits in functor $\infty$-categories are computed
    point-wise, see Corollary 5.1.2.3 in \cite{Lur09}.}, see Proposition 1.1.3.1
  in \cite{Lur17}.

\begin{proposition}
  \label{prop:pointwise} Let $\A$ be a stable $\infty$-category and $K$ a small
  $\infty$-category. Then, the $\infty$-category $\Fun[K][\A]$ of functors
  $K\to\A$ is also stable.
\end{proposition}

\section{The higher Waldhausen $\sS$-constructions}

In this section we establish a link between the higher dimensional Auslander
algebras of type $\AA$ and the higher Waldhausen $\sS$-constructions introduced
by Dyckerhoff \cite{Dyc17} and Poguntke \cite{Pog17}, and by Hesselholt and
Madsen \cite{HM15} in the case $m=2$. We begin by reminding the reader of what
is known in the classical situation, that is in the case $m=1$.

\subsection{The Waldhausen $\sS$-construction}

The following construction is due to Waldhausen \cite{Wal85}. It is the main
ingredient in the definition of the algebraic $K$-theory space of a stable
$\infty$-category\footnote{The algebraic $K$-theory space is a fundamental
  invariant which can be seen as a substantial refinement of the Grothendieck
  group.}, see for example \cite{BGT13}.

\begin{definition}
  Let $\C$ be a stable $\infty$-category and $n\geq0$. We let $\sS(\C)_n$ be the
  (stable) $\infty$-category of diagrams of the form\footnote{Such a diagram is
    precisely a functor $X\colon\Delta(1,n)\to\C$, where we endow the set
    $\Delta(1,n)$ of monotone maps $\ord 1\to\ord n$ with the natural partial
    order.}
  \begin{equation*}
    \begin{tikzcd}[column sep=tiny,row sep=tiny]
      X_{00}\rar&X_{01}\rar\dar&X_{02}\rar\dar&\cdots\rar&X_{0,n-1}\rar\dar&X_{0n}\dar\\
      &X_{11}\rar&X_{12}\rar\dar&\cdots\rar&X_{1,n-1}\rar\dar&X_{1n}\dar\\
      &&\ddots&\ddots&\vdots\dar&\vdots\dar\\
      &&&&X_{n-1,n-1}\rar&X_{n-1,n}\dar\\
      &&&&& X_{nn}
    \end{tikzcd}
  \end{equation*}
  which satisfy the following two conditions:
  \begin{itemize}
  \item For each $i\in\ord n$ the object $X_{ii}$ is a zero object of $\C$;
  \item for each $0\leq i<j<k\leq n$ the square
    \begin{equation*}
      \begin{tikzcd}
        X_{ij}\rar\dar\ar[phantom]{dr}[description]{\square}&X_{ik}\dar\\
        X_{ii}\rar&X_{jk}
      \end{tikzcd}
    \end{equation*}
    is a (homotopy) pushout diagram and a (homotopy) pullback
    diagram.\footnote{The exactness conditions imposed on the objects of
      $\sS(\C)_n$ should be thought of as a categorification of the Euler
      relations
      \eqref{eq:Euler_relations_K0An}.}.  \end{itemize}
  Size issues aside, the $\infty$-categories $\sS(\C)_n$ assemble into a
  simplicial object $\sS(\C)_\bullet$, called the \emph{Waldhausen
    $\sS$-construction of $\C$}, which takes values in the $\infty$-category of
  stable $\infty$-categories and exact functors between them.
\end{definition}

The following elementary observation can be viewed as a categorification of the
isomorphism \eqref{eq:Na1n}. Proofs can be found in Lemma 7.3 in \cite{BGT13}
and Lemma 1.2.2.4 in \cite{Lur17}. A version in the related framework of stable
derivators is proven in Theorem 4.6 in \cite{GS16} by means of a version of the
knitting algorithm.

\begin{proposition}[Waldhausen]
  \label{prop:An} Let $\C$ be a stable $\infty$-category and $n\geq1$. The
  restriction functor
  \begin{equation*}
    \functor{\sS(\C)_n}{\Fun[01\to02\to\cdots 0n][\C]},
  \end{equation*}
  which sends an object $X$ of $\sS(\C)_n $ to the diagram
  \begin{equation*}
    \begin{tikzcd}[column sep=tiny]
      X_{01}\rar&X_{02}\rar&\cdots\rar&X_{0,n-1}\rar&X_{0n},
    \end{tikzcd}
  \end{equation*}
  is an equivalence of (stable) $\infty$-categories.
\end{proposition}

The following theorem extends the foregoing proposition to arbitrary
orientations of the Dynkin diagram $\AA_n$ and can be proven using combinatorial
versions of the classical reflection functors. A proof, carried out in the
related framework of stable derivators, can be found in \cite{GS16}.

\begin{theorem}[Groth--\v{S}\v{t}ov\'\i\v{c}ek]
  \label{thm:reflection-An} Let $\C$ be a stable $\infty$-category and $n\geq1$.
  Let $S$ be a slice\footnote{As left implicit above, the Hasse quiver of
    $\Delta(1,n)$ can be identified with the Auslander--Reiten quiver of the
    quiver $\AA_n$ (with additional degenerate vertices). Slices in
    $\Delta(1,n)$ can then be defined in the usual way.} in the poset
  $\Delta(1,n)$. The restriction functor
  \begin{equation*}
    \sS(\C)_\bullet\to\Fun[\underline{S}][\C]
  \end{equation*}
  is an equivalence of (stable) $\infty$-categories.
\end{theorem}

Our aim in this section is to provide higher-dimensional versions of
\autoref{prop:An} and \autoref{thm:reflection-An} in terms of certain
higher-dimensional versions of the Waldhausen $\sS$-construction.

\subsection{The higher Waldhausen $\sS$-constructions}

Before proceeding we introduce further terminology. Let $I=\ord 1$ be the poset
$\set{0\to 1}$ and $m$ a non-negative integer. An \emph{$m$-cube} in an
$\infty$-category $\C$ is a functor $X\colon I^m\to\C$. The isomorphism
\begin{equation*}
  I^{m+1}\cong I\times I^m
\end{equation*}
together with the adjunction
\begin{equation*}
  \Fun[I\times I^m][\C]\cong\Fun[I][\Fun[I^m][\C]]
\end{equation*}
allow us to view an $(m+1)$-cube in $\C$ as morphism in the $\infty$-category
$\Fun[I^m][\C]$ of $m$-cubes in $\C$. In the
case of a stable $\infty$-category, this identification allows for an inductive
treatment of hyper-cubes as illustrated by the following definition.

\begin{definition}
  We say that a $0$-cube in a stable $\infty$-category $\A$, which is nothing
  but
  an object of $\A$, is \emph{biCartesian} if it is a zero object of $\A$.
  Inductively\footnote{The notion of a biCartesian hyper-cube can be defined
    directly as a certain (homotopy) colimit diagram, see Proposition 1.2.4.13
    and Lemma 1.2.4.5 in \cite{Lur17}.}, we say that an $(m+1)$-cube $X$ in a
  stable $\infty$-category $\A$ is \emph{biCartesian} if its cofibre (taken in
  the stable $\infty$-category $\Fun[I^m][\A]$) is a biCartesian $m$-cube in
  $\A$. For example, a $1$-cube in $\A$ is biCartesian if its underlying
  morphism is an equivalence\footnote{Equivalences are the $\infty$-categorical
    analogues of isomorphisms in ordinary category theory.} in $\A$.
\end{definition}

We are now ready to state the definition of the $m$-dimensional Waldhausen
$\sS$-construction of a stable $\infty$-category \cite{HM15,Dyc17,Pog17}.

\begin{definition}
  \label{def:S2} Let $\A$ be a stable $\infty$-category and $m$ a
  non-negative integer. For $n\geq0$ we denote by $\Sm(\A)_n$ the full subcategory of $\Fun[\Delta(m,n)][\A]$
  spanned by the diagrams $X$ satisfying the following two conditions:
  \begin{itemize}
  \item For every degenerate $m$-simplex $\sigma$ in $\Delta^n$ the object
    $X_\sigma$ is a zero object of $\A$;
  \item for each non-degenerate $(m+1)$-simplex in $\Delta^n$ consider
    the $(m+1)$-cube $q\colon I^{m+1}\to\Delta(m,n)$ given by
    $q(v)_i=\sigma_{i+v_i}$ for each $v\in I^{m+1}$ and $i\in\ord m$. Then, the
    induced $(m+1)$-cube
    \begin{equation*}
      X\circ q\colon I^{m+1}\to\A
    \end{equation*}
    is biCartesian\footnote{These exactness conditions should be thought
      of as categorifications of the higher-dimensional Euler relations
      \eqref{eq:Euler_relations_Am_1}. For example, for $m=2$ the cubes which
      are required to be biCartesian are those of the form
  $$
  \begin{tikzcd}[column sep=tiny, row sep=tiny,ampersand replacement=\&]
    X_{ijk}\ar{rr}\drar\ar{dd}\&\&X_{ijl}\drar\ar{dd}\\
    \&X_{jjk}\ar[crossing over]{rr}\&\&X_{jjl}\ar{dd}\\
    X_{ikk}\ar{rr}\drar\&\&X_{ikl}\drar\\
    \&X_{jkk}\ar{rr}\ar[crossing over, leftarrow]{uu}\&\&X_{jkl}
  \end{tikzcd}
    $$
    where $0\leq i<j<k<l\leq n$; these cubes are $\infty$-categorical analogues
    of 4-term exact sequences.}.
\end{itemize}
Size issues aside, the (stable) $\infty$-categories $\Sm(\A)_n$, $n\geq 0$, assemble into a simplicial object $\Sm(\A)_\bullet$, called the \emph{$m$-dimensional Waldhausen
  $\sS$-construction of $\A$}, which takes values in the $\infty$-category of
stable $\infty$-categories and exact functors between them.
\end{definition}

\begin{remark}
  The second exactness conditions appearing in the definition of the
  $m$-dimensional Waldhausen $\sS$-construction are equivalent to the condition
  that the $(m+1)$-cubes of the form $X\circ q$ are biCartesian,
  where $$q(v)_i=\sigma_i+v_i,\quad i\in\ord m$$ and $\sigma$ ranges over all
  non-degenerate $m$-simplices in $\Delta^n$ with $\sigma_m<n$. For example, for
  $m=2$ these are the cubes
    $$
    \scalebox{0.85}{
      \begin{tikzcd}[column sep=0.25em, row sep=tiny,ampersand replacement=\&]
        X_{ijk}\ar{rr}\drar\ar{dd}\&\&X_{ij,k+1}\drar\ar{dd}\\
        \&X_{i+1,j,k}\ar[crossing over]{rr}\&\&X_{i+1,j,k+1}\ar{dd}\\
        X_{i,j+1,k}\ar{rr}\drar\&\&X_{i,j+1,k+1}\drar\\
        \&X_{i+1,j+1,k}\ar{rr}\ar[crossing over,
        leftarrow]{uu}\&\&X_{i+1,j+1,k+1}
      \end{tikzcd}
    }
    $$
    where $0\leq i<j<k<n$. These conditions should be thought of as a
    categorification of the higher Auslander--Reiten relations
    \eqref{eq:Euler_relations_K0An}.
\end{remark}

\begin{example}
  Let $\A$ be a stable $\infty$-category, $m=2$, and $n=4$. After discarding
  redundant information pertaining additional zero objects, an object of
  $\sS^{\langle 2\rangle}(\A)_4$ can be identified with a diagram of the form
  \begin{equation}
    \label{eq:S24}
    \begin{tikzpicture}
      \node (012) at (0,0,0) {$X_{012}$}; \node (013) at (2,0,0) {$X_{013}$};
      \node (014) at (4,0,0) {$X_{014}$};

      \node (022) at (0,-2,0) {$0$}; \node (023) at (2,-2,0) {$X_{023}$}; \node
      (024) at (4,-2,0) {$X_{024}$};

      \node (033) at (2,-4,0) {$0$}; \node (034) at (4,-4,0) {$X_{034}$};

      \node (112) at (0,0,2) {$0$}; \node (113) at (2,0,2) {$0$}; \node (114) at
      (4,0,2) {$0$};

      \node (122) at (0,-2,2) {$0$}; \node (123) at (2,-2,2) {$X_{123}$}; \node
      (124) at (4,-2,2) {$X_{124}$};

      \node (133) at (2,-4,2) {$0$}; \node (134) at (4,-4,2) {$X_{134}$};

      \node (223) at (2,-2,4) {$0$}; \node (224) at (4,-2,4) {$0$}; \node (233)
      at (2,-4,4) {$0$}; \node (234) at (4,-4,4) {$X_{234}$};
        
      \path[commutative diagrams/.cd, every arrow]%
      (012) edge (013)%
      (013) edge (014)%
      (012) edge (022)%
      (013) edge (023)%
      (014) edge (024)%
      (022) edge (023)%
      (023) edge (024)%
      (023) edge (033)%
      (024) edge (034)%
      (033) edge (034)%
      (012) edge (112)%
      (013) edge (113)%
      (014) edge (114)%
      (022) edge (122)%
      (023) edge (123)%
      (024) edge (124)%
      (033) edge (133)%
      (034) edge (134)%
      ; \path[commutative diagrams/.cd, every arrow, every label]%
      (112) edge[commutative diagrams/crossing over] (113)%
      (113) edge[commutative diagrams/crossing over] (114)%
      (112) edge (122)%
      (113) edge[commutative diagrams/crossing over] (123)%
      (114) edge[commutative diagrams/crossing over] (124)%
      (122) edge (123)%
      (123) edge[commutative diagrams/crossing over] (124)%
      (123) edge (133)%
      (124) edge[commutative diagrams/crossing over] (134)%
      (133) edge (134)%
      (123) edge (223)%
      (124) edge (224)%
      (133) edge (233)%
      (134) edge (234)%
      ; \path[commutative diagrams/.cd, every arrow]%
      (223) edge[commutative diagrams/crossing over] (224)%
      (223) edge (233)%
      (224) edge[commutative diagrams/crossing over] (234)%
      (233) edge (234)%
      ;
    \end{tikzpicture}
  \end{equation}  
  in which all `unit cubes' are biCartesian. A standard application of the
  theory of Kan extensions in $\infty$-categories shows that such a diagram is
  determined by its restriction to the (coherently commutative) diagram
  \begin{equation*}
    \begin{tikzcd}[row sep=large]
      X_{012}\rar\dar&X_{013}\rar\dar&X_{014}\dar\\
      0\rar&X_{023}\rar\dar&X_{024}\dar\\
      &0\rar&X_{034}
    \end{tikzcd}
  \end{equation*}
  The latter diagram can be thought of as a (coherent) representation of the
  Auslander algebra $\AA_3^{(2)}$ of the quiver $\AA_3$ in the stable
  $\infty$-category $\A$. It is also worth noting that diagram \eqref{eq:S24}
  agrees with the `higher-dimensional Auslander--Reiten quiver' of the unique
  clus\-ter-til\-ting subcategory associated with $\AA_3^{(2)}$, see Section 6
  in \cite{Iya11} or \cite{OT12}.
\end{example}

The following observation relates the $m$-dimensional Waldhausen
$\sS$-construction to the $m$-dimensional Auslander algebras of type $\AA$; it
can be viewed as a higher-dimensional analogue of \autoref{prop:An}. It can be
proven either constructively, by means of a higher-dimensional version of the
knitting algorithm\footnote{This approach is analogous to the proof of Theorem
  4.6 in \cite{GS16}. The aforementioned knitting algorithm can be deduced, for
  example, from the proof of Theorem 5.27 in \cite{IO11}.}, or by means of the
point -wise formulas for $\infty$-categorical Kan extensions\footnote{This
  approach is analogous to the proof of \autoref{prop:An} given for example in
  Lemma 1.2.2.4 in \cite{Lur17}.}. A proof can be found in Proposition 2.10 in
\cite{DJW19b}.

\begin{proposition}[Dyckerhoff--J--Walde]
  \label{prop:Ank} Let $\A$ be a stable $\infty$-category and $n\geq m\geq1$.
  The restriction functor
  \begin{equation*}
    \functor{\Sm(\A)_m}{\Fun[P(m,n)][\A][][*]}
  \end{equation*}
  is an equivalence of (stable) $\infty$-categories, where
  \begin{equation*}
    P(m,n):=\setP{\sigma\in\Delta(m,n)}{\sigma_0=0}
  \end{equation*}
  and $\Fun[P(m,n)][\A][][*]$ is the full subcategory of $\Fun[P(m,n)][\A]$
  spanned by those functors which send degenerate $m$-simplices in $\Delta^n$
  to zero objects in $\A$.
\end{proposition}

\begin{remark}
  After discarding further redundant information, the poset $P(m,n)$ appearing
  in in \autoref{prop:Ank} can be replaced by a smaller poset which models
  perfectly the quiver with relations of the higher Auslander algebra
  $\AA_{n-m+1}^{(m)}$. This fact is perhaps more transparent from the
  description of Iyama's higher Auslander algebras of type $\AA$ given in
  Section 2 in \cite{JK16}.
\end{remark}

As explained in \cite{IO11}, the following theorem can be regarded as a
higher-dimensional version of \autoref{thm:reflection-An}. A proof can be found
in Theorem 2.41 in \cite{DJW19b}.

\begin{theorem}[Dyckerhoff--J--Walde]
  \label{thm:slices} Let $\A$ be a stable $\infty$-category and $n\geq m\geq1$.
  Let $S$ be a slice\footnote{The notion of a `slice' is analogous to that
    introduced in \cite{IO11}.} in the poset $\Delta(m,n)$. The restriction
  functor
  \begin{equation*}
    \Sm(\A)_n\to\Fun[\underline{S}][\A][][*]
  \end{equation*}
  is an equivalence of (stable) $\infty$-categories, where $\underline{S}$
  denotes the convex hull of $S$ in the poset $\Delta(m,n)$.
\end{theorem}

\begin{remark}
  The proof of \autoref{thm:slices} relies on combinatorial versions of the
  derived equivalences induced by higher-dimensional reflection
  functors\footnote{Combinatorial versions of classical reflection functors are
    investigated in \cite{GS16} in the case of quivers of type $\AA$ and in
    \cite{GS16a} for general trees in the related framework of stable
    derivators. A further generalisation, corresponding to the generalised
    reflection functors of Ladkani \cite{Lad08}, was obtained in \cite{DJW19a}.}
  in the sense of \cite{IO11}; these reflection functors rely on the operation
  of slice mutation also introduced in \cite{IO11}. In more detail, if $S$ and
  $S'$ are slices in $\Delta(m,n)$ which are mutation of each other, then there
  exists a slightly larger poset $\underline{S\diamond S}'$ containing both $S$
  and $S'$ as well as a distinguished $(m+1)$-cube. This larger poset allows us
  to realise the aforementioned reflection functors by means of equivalences of
  stable $\infty$-categories
  $$
    \begin{tikzcd}[ampersand replacement=\&,column sep=small, row sep=small]
      \Fun[\underline{S}][\A][][*]\&\Fun[\underline{S\diamond S}'][\A][ex][*]\lar[swap]{\simeq}\rar{\simeq}\&\Fun[\underline{S}'][\A][][*]
    \end{tikzcd}
    $$
    induced by the restriction functors, where the objects of the stable
    $\infty$-category $\Fun[\underline{S\diamond S}'][\A][ex][*]$ satisfy the additional
    requirement that their restriction along the distinguished $(m+1)$-cube in
    $\underline{S\diamond S}'$ is a biCartesian $(m+1)$-cube in $\A$, see Figure
    \ref{fig:slices}. The proof of \autoref{thm:slices} is obtained by
    combining these equivalences with \autoref{prop:Ank} and the transitivity of
    the mutation operation on slices proven in \cite{IO11}. See \cite{DJW19b} for
    details.
    \begin{figure*}
      \begin{center}
      \scalebox{0.7}{        
      \begin{tikzpicture}
        \node (013) at (2,0,0) {$X_{013}$};%
        \node (014) at (4,0,0) {$X_{014}$};%
        \node (023) at (2,-2,0) {$X_{023}$};%
        \node (024) at (4,-2,0) {$X_{024}$};%
        \node[rotate=180] (025) at (5,-2,0) {$\mapsto$}; \node (026) at (8,-2,0)
        {$\phantom{\mapsto}$}; \node (033) at (2,-4,0) {$0$};%
        \node (034) at (4,-4,0) {$X_{034}$};%
        \path[commutative diagrams/.cd, every arrow]%
        (013) edge (014)%
        (013) edge (023)%
        (014) edge (024)%
        (023) edge (024)%
        (023) edge (033)%
        (024) edge (034)%
        (033) edge (034)%
        ;

        \node (022) at (0,-2,0) {$0$};%

        \path[commutative diagrams/.cd, every arrow]%
        (022) edge (023)
        ;

        \node (122) at (0,-2,2) {$\phantom{0}$}; \node (012) at (0,0,0)
        {$X_{012}$}; \path[commutative diagrams/.cd, every arrow]%
        (012) edge (013)%
        (012) edge (022)%
        ; \node at (2,1.0,0) {$\Fun[\underline{S}][\A][][*]$};
        \begin{scope}[xshift=20em]
        \node (013) at (2,0,0) {$X_{013}$};%
        \node (014) at (4,0,0) {$X_{014}$};%
        \node (023) at (2,-2,0) {$X_{023}$};%
        \node (024) at (4,-2,0) {$X_{024}$};%

        \node (033) at (2,-4,0) {$0$};%
        \node (034) at (4,-4,0) {$X_{034}$};%
        \node (123) at (2,-2,2) {$X_{123}$};%
        \node (113) at (2,0,2) {$0$};%
        \path[commutative diagrams/.cd, every arrow]%
        (013) edge (014)%
        (013) edge (023)%
        (014) edge (024)%
        (023) edge (024)%
        (023) edge (033)%
        (024) edge (034)%
        (033) edge (034)%
        (023) edge (123)%
        (013) edge (113)%
        (113) edge[commutative diagrams/crossing over] (123)%
        ;

        \node (112) at (0,0,2) {$0$}; \node (122) at (0,-2,2) {$0$}; \node (022)
        at (0,-2,0) {$0$};%

        \path[commutative diagrams/.cd, every arrow]%
        (022) edge (023) (122) edge (123)%
        (112) edge[commutative diagrams/crossing over] (113)%
        (022) edge (122)%
        (112) edge (122)%
        ;

        \path[commutative diagrams/.cd, every arrow]%
        (113) edge[commutative diagrams/crossing over] (123)%
        ;
        
        \node (012) at (0,0,0) {$X_{012}$}; \node (025) at (5.25,-2,0) {$\mapsto$};
        \path[commutative diagrams/.cd, every arrow]%
        (023) edge (123)%
        (122) edge (123)%
        (012) edge (013)%
        (012) edge (112)%
        (012) edge (022)%
        (112) edge[commutative diagrams/crossing over] (113)%
        ; \node at (2,1.0,0) {$\Fun[\underline{S\diamond S}'][\A][ex][*]$};
        \begin{scope}[xshift=15em]
        \node (013) at (2,0,0) {$X_{013}$};%
        \node (014) at (4,0,0) {$X_{014}$};%
        \node (023) at (2,-2,0) {$X_{023}$};%
        \node (024) at (4,-2,0) {$X_{024}$};%

        \node (033) at (2,-4,0) {$0$};%
        \node (034) at (4,-4,0) {$X_{034}$};%
        \node (123) at (2,-2,2) {$X_{123}$};%
        \node (113) at (2,0,2) {$0$};%
        \path[commutative diagrams/.cd, every arrow]%
        (013) edge (014)%
        (013) edge (023)%
        (014) edge (024)%
        (023) edge (024)%
        (023) edge (033)%
        (024) edge (034)%
        (033) edge (034)%
        (023) edge (123)%
        (013) edge (113)%
        (113) edge[commutative diagrams/crossing over] (123)%
        ;



        \node (012) at (0,0,0) {$\phantom{X_{012}}$};
        \node at (2,1.0,0) {$\Fun[\underline{S}'][\A][][*]$};
      \end{scope}
    \end{scope}        
      \end{tikzpicture}
    }
    \caption{An example of slice mutation.}
    \label{fig:slices}
  \end{center}
\end{figure*}
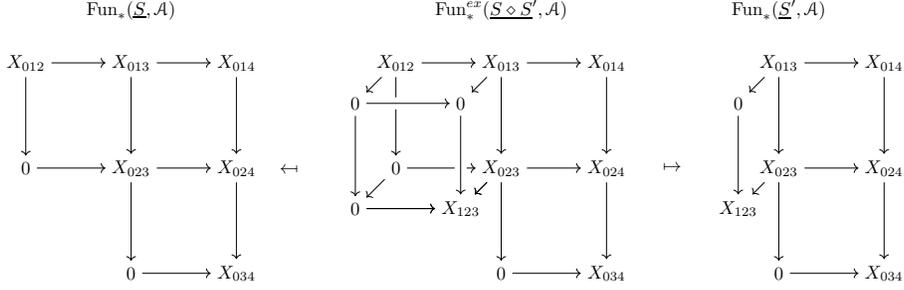
\end{remark}

\subsection{Recollements}

We make the following general observation, which makes the inductive nature of
the $m$-dimensional Waldhausen $\sS$-construction of $\A$ readily apparent.
Indeed, \autoref{prop:ladder} is reminiscent of Iyama's inductive description of
the higher Auslander algebras of type $\AA$ by means of cones of translation
quivers, see Section 6 in \cite{Iya11}. A proof can be found in Proposition 2.51
in \cite{DJW19b}.

\begin{proposition}[Dyckerhoff--J--Walde]
  \label{prop:ladder} Let $\A$ be a stable $\infty$-category and $n\geq m\geq
  1$ integers. For each $i\in\ord n$ the functor
  $s_i\colon\Sm(\A)_n\to\Sm(\A)_{n+1}$ is part of a recollement of stable
  $\infty$-categories
  \begin{equation*}
    \recollementnoph{\Sm(\A)_n}{\Sm(\A)_{n+1}}{\sS^{\langle m-1\rangle}(\A)_n.}{d_i}{s_i}{d_{i+1}}
  \end{equation*}
  In particular, the sequence of adjunctions
  \begin{equation*}
    d_0\dashv s_0\dashv d_1\dashv s_1\dashv\cdots\dashv s_n\dashv d_{n+1}
  \end{equation*}
  is part of a ladder of recollements in the sense of \cite{BBD82,AHKLY17}.
\end{proposition}

\section{Acknowledgements}

The author would like to acknowledge the numerous conversations with Julian
K{\"u}lshammer concerning the combinatorial treatment of the higher-dimensional
Auslander algebras of type $\AA$ which helped shape the author's viewpoint on
the subject. The author also thanks Marius Thaule and Tashi Walde for comments
on a preliminary version of this note. Versions of \autoref{prop:Ank} and
\autoref{prop:ladder} have been obtained by Beckert in his PhD thesis
\cite{Bec18} in the related framework of stable derivators. The author was
funded by the Deutsche Forschungs-gemeinschaft (DFG, German Research Foundation)
under Germany's Excellence Strategy - GZ 2047/1, Projekt-ID 390685813.

\providecommand{\bysame}{\leavevmode\hbox to3em{\hrulefill}\thinspace}
\providecommand{\MR}{\relax\ifhmode\unskip\space\fi MR }
\providecommand{\MRhref}[2]{%
  \href{http://www.ams.org/mathscinet-getitem?mr=#1}{#2}
}
\providecommand{\href}[2]{#2}

\end{document}